\documentclass[a4paper,12pt]{article}
\usepackage{comment}
\usepackage{cite}
\usepackage{amsmath}
\usepackage{amssymb}
\usepackage{amsfonts}
\usepackage[T1]{fontenc}
\usepackage[utf8]{inputenc}
\usepackage{graphicx}
\usepackage{fancyhdr}
\usepackage{float}
\usepackage{xcolor}
\usepackage{authblk}
\usepackage{mathrsfs}
\usepackage{empheq}
\usepackage{url}

\usepackage{geometry}
\begin{document}

\title{New series representations for any positive power of $\pi$ from a relation involving trigonometric functions}

\author[$\dagger$]{Jean-Christophe {\sc Pain}\\
\small
CEA, DAM, DIF, F-91297 Arpajon, France\\
Universit\'e Paris-Saclay, CEA, Laboratoire Mati\`ere en Conditions Extr\^emes,\\ 
91680 Bruy\`eres-le-Ch\^atel, France
}

\maketitle

\begin{abstract}
In previous works, we presented series representations for $\pi^3$ and $\pi^5$, in which the prefactor depends only on the golden ratio appears. In this article, we derive a general relation involving trigonometric functions and an infinite series. Such an identity is likely to provide many series representations for any positive power of $\pi$, among them the above mentioned representations for $\pi^3$ and $\pi^5$.
\end{abstract}

\section{Introduction}

There are very few known expressions for $\pi^m$ with $m=2k+1$ and $k\in\mathbb{N}^*$ \cite{Borwein1987,Pain2022a,Pain2022b}. We recently published series representations for $\pi^3$ \cite{Pain2022c} and $\pi^5$ \cite{Pain2022d}, consisting of a prefactor depending only on the golden ratio multiplied by an infinite sum. The latter results were both obtained as particular cases of identities involving trigonometric functions. In the present article, we show that such identities are in fact particular cases of a more general equation, providing in that way a family of series for any positive power of $\pi$. The general formula for $\pi^{k+2}$ with $k\ge 0$ involves trigonometric functions, multinomial coefficients and an infinite series. It is derived in section \ref{sec2}, and the special cases corresponding to $k=1$ and $k=2$ are obtained in section \ref{sec3}.

\section{General trigonometric relation providing series representations for $\pi^{k+2}$, $k\ge 0$}\label{sec2}

We start from (see Ref. \cite{Borwein1987}, 13.a, p. 382):
\begin{equation}
\pi\left[\cot(\pi x)-\cot(\pi a)\right]=\sum_{-\infty}^{\infty}\frac{a-x}{(x-n)(a-n)}.
\end{equation}
Differentiating with respect to $x$ yields
\begin{equation}\label{comp}
\frac{\pi^2}{\sin^2(\pi x)}=\sum_{n=-\infty}^{\infty}\frac{1}{(x-n)^2}.
\end{equation}
Differentiating once more with respect to $x$ gives ($\mathrm{cosec}(y)=1/\sin(y)$)
\begin{equation}\label{pitrois}
\pi^3\cot(\pi x)\mathrm{cosec}^2(\pi x)=\sum_{n=-\infty}^{\infty}\frac{1}{(x-n)^3},
\end{equation}
which is exactly formula 13.b, p. 382 in Ref. \cite{Borwein1987}. Taking the $k^{th}$ derivative of both sides of Eq. (\ref{comp}) leads to

\begin{equation}\label{base}
\pi^2\frac{\partial^k}{\partial x^k}\left[\frac{1}{\sin^2(\pi x)}\right]=\sum_{n=-\infty}^{\infty}\frac{(-1)^k(k+1)!}{(x-n)^{k+2}}.
\end{equation}
The Fa\`a di Bruno formula \cite{Arbogast1800,Faa1855,Faa1857,Cesaro1885,Comtet1974,Roman1980} tells us that, if $f(x)=g(h(x))$, then, using the notation
\begin{equation}
\frac{\partial^k}{\partial x^k}f(x)=f^{(k)}(x),
\end{equation}
we have
\begin{equation}\label{faa}
f^{(k)}(x)=\sum_{\left\{m_i\right\}}\frac{k!}{m_1!m_2!\cdots m_k!}g^{(m_1+m_2+\cdots m_k)}(h(x))\prod_{j=1}^k\left[\frac{h^{(j)}(x)}{j!}\right]^{m_j}
\end{equation}
with
\begin{equation}
\sum_{j=1}^kjm_j=k.
\end{equation}
In the present case, $g(x)=1/x$, $h(x)=\sin^ 2(\pi x)$ and
\begin{equation}
g^{(p)}(x)=\frac{(-1)^pp!}{x^{p+1}}
\end{equation}
as well as
\begin{equation}
h^{(p)}(x)=-2^{p-1}~\pi^p~\cos\left(2\pi x+p\frac{\pi}{2}\right).
\end{equation}
Therefore, Eq. (\ref{base}) becomes
\begin{equation}
\pi^{2+k}\sum_{\substack{\sum_{j=1}^kjm_j=k}}\frac{\left(\sum_{j=1}^km_j\right)!~2^k\prod_{j=1}^k\left[\cos\left(2\pi x+j\frac{\pi}{2}\right)\right]^ {m_j}}{2^{\sum_{j=1}^km_j}\left(\prod_{j=1}^km_j!\right)\left[\sin^2(\pi x)\right]^ {\sum_{j=1}^km_j+1}\prod_{j=1}^k(j!)^{m_j}}=\sum_{n=-\infty}^{\infty}\frac{(-1)^k(k+1)}{(x-n)^{k+2}},
\end{equation}
which can be put in the form
\begin{empheq}[box=\fbox]{align}
\pi^{2+k}=\frac{(-1)^k(k+1)}{2^k\mathscr{A}_k(x)}\sum_{n=-\infty}^{\infty}\frac{1}{(x-n)^{k+2}},
\end{empheq}
where we have set
\begin{empheq}[box=\fbox]{align}
\mathscr{A}_k(x)=\sum_{\left\{m_i\right\}}\frac{\mathscr{S}_k!}{2^{\mathscr{S}_k}\left(\prod_{j=1}^km_j!\right)\left[\sin^2(\pi x)\right]^ {\mathscr{S}_k+1}}\times\prod_{j=1}^k\frac{\left[\cos\left(2\pi x+j\frac{\pi}{2}\right)\right]^{m_j}}{(j!)^{m_j}},
\end{empheq}
and
\begin{empheq}[box=\fbox]{align}
\mathscr{S}_k=\sum_{j=1}^km_j.
\end{empheq}
Choosing a particular value of $x$ (for instance $x=1/4$), one obtains an infinite series representation of $\pi^{k+2}$, whatever the positive integer value of $k$. In particular, choosing $x=1/5$ or $1/10$ \cite{Pain2022c,Pain2022d}, it is possible to obtain infinite series in which the prefactor involves the golden ratio.

\section{Examples}\label{sec3}

For $k=1$, we have only one set of $\left\{m_i\right\}$ corresponding to the singleton $m_1=1$. In that case $\mathscr{S}_1=1$ and we recover easily Eq. (\ref{pitrois}).
In the case $k=2$, we have two possibilities (two possible pairs): $(m_1,m_2)=(2,0)$ (corresponding to $\mathscr{S}_2=2$) and $(m_1,m_2)=(0,1)$ (corresponding to $\mathscr{S}_2=1$). The former yields
\begin{equation}
\mathscr{T}_1=\pi^4\frac{\cos^2\left(2\pi x+\frac{\pi}{2}\right)}{\sin^6(\pi x)}=\pi^4\frac{\sin^2(2\pi x)}{\sin^6(\pi x)}=\pi^4\frac{4\sin^2(\pi x)\left(1-\sin^2(\pi x)\right)}{\sin^6(\pi x)}
\end{equation}
and the latter
\begin{equation}
\mathscr{T}_2=-\pi^4\frac{\cos\left(2\pi x+\pi\right)}{\sin^4(\pi x)}=\pi^ 4\left(-\frac{1}{\sin^4(\pi x)}+\frac{2}{\sin^4(\pi x)}\right),
\end{equation}
and thus
\begin{equation}
\mathscr{T}_1+\mathscr{T}_2=3\sum_{n=-\infty}^{\infty}\frac{1}{(x-n)^4}.
\end{equation}
This gives
\begin{equation}
\pi^4\left[\frac{4}{\sin^4(\pi x)}-\frac{4}{\sin^2(\pi x)}-\frac{1}{\sin^4(\pi x)}+\frac{2}{\sin^2(\pi x)}\right]=3\sum_{n=-\infty}^{\infty}\frac{1}{(x-n)^4},
\end{equation}
yielding
\begin{equation}
\pi^4\left[\mathrm{cosec}^4(\pi x)-\frac{2}{3}\mathrm{cosec}^2(\pi x)\right]=\sum_{n=-\infty}^{\infty}\frac{1}{(x-n)^4},
\end{equation}
which is exactly Eq. 13.c p. 382 of Ref. \cite{Borwein1987} and can be easily obtained by deriving Eq. (\ref{pitrois}) with respect to $x$.

\section{Conclusion}

We recently published series representations for $\pi^3$ and $\pi^5$, in which the golden ratio appears in the prefactor. In the present work, we obtained a general relation involving trigonometric functions and an infinite series. The latter identity is likely to provide families of series representations for any power of $\pi$. The above mentioned previously-obtained representations for $\pi^3$ and $\pi^5$ are particular cases of the present formalism. The new relation, obtained as a multiple derivative of an identity expressing the difference of two cotangent functions as an infinite series, is obtained using the Fa\`a di Bruno formula for the multiple derivative of a composite function. The resulting expression may also provide relations between powers of $\pi$ and other mathematical constants.

\section{Appendix: on the Fa\`a di Bruno formula}

\subsection{Proof by recurrence}

Let us prove by recurrence the expression
\begin{equation}\label{tbd}
\left[g(h(x))\right]^{(k)}=\sum_{\substack{\sum_{j=1}^kjp_j=k}}\mathscr{C}(k,m_1,m_2,\cdots,m_k)~g^{\sum_{j=1}^km_j}(h(x))~\prod_{j=1}^k\left[h^{(j)}(x)\right]^{m_j},
\end{equation}
where $\mathscr{C}$ is a coefficient which will be precised later. For $k=1$, $m_1=1$ and $\left[g(h(x))\right]'=g'(h(x))\times f'(x)$, which means that the expression (\ref{tbd}) is true. One has

\begin{eqnarray}
\left[g(h(x))\right]^{(k+1)}&=&\sum_{\left\{m_i\right\}}\mathscr{C}(k,m_1,m_2,\cdots,m_k)~g^{1+\sum_{j=1}^km_j}(h(x))\times h'(x)~\prod_{j=1}^k\left[h^{(j)}(x)\right]^{m_j}\nonumber\\
& &+\sum_{\left\{m_i\right\}}\mathscr{C}(k,m_1,m_2,\cdots,m_k)~g^{\sum_{j=1}^km_j}(h(x))~\left\{\prod_{j=1}^k\left[h^{(j)}(x)\right]^{m_j}\right\}'.
\end{eqnarray}
The first term in the right-hand side reads
\begin{equation}\label{first}
\sum_{\substack{\sum_{j=1}^kjn_j=k+1}}\mathscr{C}(k,m_1,m_2,\cdots,m_k)~g^{\sum_{j=1}^kn_j}(h(x))~\prod_{j=1}^k\left[h^{(j)}(x)\right]^{n_j},
\end{equation}
where $n_1=m_1+1$ and $n_i=m_i$ for $i\ge 2$. The quantity (\ref{first}) is of the required form, bringing a contribution $\mathscr{C}(k,m_1,m_2,\cdots,m_k)$ to the term $\mathscr{C}(k+1,m_1+1,m_2,\cdots,m_k,0)$. Since we have
\begin{equation}
\left\{\prod_{j=1}^k\left[h^{(j)}(x)\right]^{m_j}\right\}'=\left[\sum_{i=1}^km_i\frac{h^{(i+1)}(x)}{h^{(i)}(x)}\right]\prod_{j=1}^k\left[h^{(j)}(x)\right]^{m_j},
\end{equation}
the second sum is made of terms proportional to
\begin{equation}
g^{\sum_{j=1}^km_j}(h(x))\times \prod_{j=1}^k\left(h^{(j)}(x)\right)^{p_j},
\end{equation}
with $p_i=m_i-1$, $p_{i+1}=m_{i+1}+1$ and $p_j=m_j$ for all $j$ different from $i$ or $i+1$. Thus
\begin{equation}
\sum_{j=1}^km_j=\sum_{j=1}^{k+1}p_j
\end{equation}
and
\begin{equation}
\sum_{j=1}^{k+1}jp_j=\sum_{j=1}^{k}jp_j+(i+1)-i=k+1.
\end{equation}
This shows that the recurrence is true. More precisely, we found that $\mathscr{C}(k+1,m_1,m_2,\cdots,m_k,0)$ is sum of $\mathscr{C}(k,m_1-1,m_2,\cdots,m_k)$ and terms of the form $(m_i+1)\mathscr{C}(k+1,m_1,m_2,\cdots,m_i+1,m_{i+1}-1,\cdots,m_k)$. We have also $\mathscr{C}(k,m_1,m_2,\cdots,m_k)=1$. We now have to prove that
\begin{equation}\label{coeff}
\mathscr{C}(k,m_1,m_2,\cdots,m_k)=\frac{k!}{\prod_{j=1}^k\left(m_j!\times (j!)^{m_j}\right)}
\end{equation}
for a suite such that
\begin{equation}
\sum_{j=1}^kjm_j=k.
\end{equation}
In other words, one needs to show that
\begin{equation}
\mathscr{C}(k+1,n_1,n_2,\cdots,n_k,0)=\frac{(k+1)!}{\prod_{j=1}^k\left(n_j!\times (j!)^{n_j}\right)}.
\end{equation}
We have
\begin{eqnarray}
\mathscr{C}(k+1,n_1,n_2,\cdots,n_k,0)&=&\mathscr{C}(k,n_1-1,n_2,\cdots,n_k)\nonumber\\
& &+(n_1+1)~\mathscr{C}(k,n_1+1,n_2-1,\cdots,n_k)\nonumber\\
& &+(n_2+1)~\mathscr{C}(k,n_1,n_2+1,n_3-1,\cdots,n_k)\nonumber\\
& &+\cdots\nonumber\\
& &+(n_{k-1}+1)~\mathscr{C}(k,n_1,n_2+1,\cdots,n_{k-1}+1,n_k-1).
\end{eqnarray}
Since one has
\begin{equation}
n_1+2n_2+3n_3+\cdots +kn_k=k+1,
\end{equation}
one gets
\begin{eqnarray}
\frac{(k+1)!}{\prod_{j=1}^k\left(n_j!\times (j!)^{n_j}\right)}&=&\frac{k!\times n_1\times 1!}{\prod_{j=1}^k\left(n_j!\times (j!)^{n_j}\right)}\nonumber\\
& &+\frac{k!\times (n_1+1)\times n_2\times 2!}{(n_1+1)\times 1!\times\prod_{j=1}^k\left(n_j!\times (j!)^{n_j}\right)}\nonumber\\
& &+\frac{k!\times (n_2+1)\times n_3\times 3!}{(n_2+1)\times 2!\times\prod_{j=1}^k\left(n_j!\times (j!)^{n_j}\right)}\nonumber\\
& &+\cdots\nonumber\\
& &+\frac{k!\times (n_{k-1}+1)\times n_k\times k!}{(n_{k-1}+1)\times (k-1)!\times\prod_{j=1}^k\left(n_j!\times (j!)^{n_j}\right)},
\end{eqnarray}
which completes the proof, justifying the expression (\ref{coeff}).

\subsection{Alternative formula for the multiple derivative of the inverse of a function} 

In the case of the $p^{th}$ derivative of the inverse of a function, the Fa\`a di Bruno formula can be put in the form
\begin{equation}
\left(\frac{1}{g}\right)^{(p)}(x)=\frac{p!}{g^{p+1}(x)}\sum_{\left\{m_i\right\}}\frac{(-1)^{p-p_0}(p-p_0)!}{\prod_{j=1}^p(j!)^{m_j}m_j!}\prod_{j=0}^p\left[g^{(j)}(x)\right]^{m_j},
\end{equation}
with $p_0=m_2+2m_3+\cdots +(p-1)m_p$ and $p=m_0+m_1+\cdots +m_p=m_1+2m_2+\cdots +pm_p$.

\subsection{Alternative expressions in terms of Bell polynomials}

Formula (\ref{faa}) can also be written as
\begin{equation}
\left[g(h(x))\right]^{(k)}=\sum_{\ell=0}^kg^{(\ell)}(h(x))~B_{k,\ell}\left(h'(x),h''(x),\cdots, h^{(k-\ell+1)}(x)\right),
\end{equation}
where the partial or ordinary Bell polynomials $B_{n,\ell}$ \cite{Bell1934} are a triangular array of polynomials given by
\begin{equation}
B_{k,\ell}\left(x_1,x_2,\cdots,x_{k-\ell+1}\right)=\sum_{\left\{j_i\right\}}\frac{k!}{j_1!j_2!\cdots j_{k-\ell+1}!}\left(\frac{x_1}{1!}\right)^{j_1}\left(\frac{x_2}{2!}\right)^{j_2}\cdots \left(\frac{x_{k-\ell+1}}{(k-\ell+1)!}\right)^{j_{k-\ell+1}},
\end{equation}
with $j_1+j_2+\cdots +j_{k-\ell+1}=\ell$ and $j_1+2j_2+\cdots +(k-\ell+1)j_{k-\ell+1}=k$. One has
\begin{equation}
B_k\left(x_1,x_2,\cdots, x_k\right)=\sum_{\ell=1}^kB_{k,\ell}(x_1,x_2,\cdots,x_{k-\ell+1}),
\end{equation}
where $B_k$ are the complete exponential Bell polynomials satisfying
\begin{equation}
\exp\left(\sum_{i=1}^{\infty}\frac{a_i}{i!}x^i\right)=\sum_{n=0}^{\infty}\frac{B_n(a_1,a_2,\cdots,a_n)}{n!}x^n.
\end{equation}

\end{document}